\documentclass[reqno]{amsart}

\usepackage{amsmath,amssymb,amsthm,amsfonts,mathrsfs}
\usepackage{fullpage}
\usepackage{xcolor}
\usepackage{graphicx}
\usepackage{listings}

\definecolor{codegreen}{rgb}{0,0.6,0}
\definecolor{codegray}{rgb}{0.5,0.5,0.5}
\definecolor{codepurple}{rgb}{0.58,0,0.82}
\definecolor{backcolour}{rgb}{0.95,0.95,0.92}

\lstdefinestyle{mystyle}{
  backgroundcolor=\color{backcolour},   commentstyle=\color{codegreen},
  keywordstyle=\color{magenta},
  numberstyle=\tiny\color{codegray},
  stringstyle=\color{codepurple},
  basicstyle=\ttfamily\footnotesize,
  breakatwhitespace=false,         
  breaklines=true,                 
  captionpos=b,                    
  keepspaces=true,                 
  numbers=left,                    
  numbersep=5pt,                  
  showspaces=false,                
  showstringspaces=false,
  showtabs=false,                  
  tabsize=2
}

\usepackage{listings}
\usepackage[colorlinks=true,
linkcolor=blue,
anchorcolor=blue,
citecolor=red
]{hyperref}
\allowdisplaybreaks % allow equations break between two pages
\newtheorem{theorem}{Theorem}[section]
\newtheorem{lemma}[theorem]{Lemma}

\newtheorem*{conjecture*}{Conjecture}
\theoremstyle{definition}

\theoremstyle{remark}

\newtheorem*{remark*}{remark}
\renewcommand{\leq}{\leqslant}

\usepackage{colonequals}

\author{Runbo Li}
\address{International Curriculum Center, The High School Affiliated to Renmin University of China, Beijing, China}
\email{runbo.li.carey@gmail.com}

\makeatletter
\@namedef{subjclassname@2020}{\textup{}2020 Mathematics Subject Classification}
\makeatother

\title[]{On the exceptional set in the $abc$ conjecture}
\subjclass[2020]{11D45, 11D75}
\keywords{exceptional set, Diophantine equation, combinatorial optimization}

\begin{document}
	
\begin{abstract}
The $abc$ conjecture states that there are only finitely many triples of coprime positive integers $(a,b,c)$ such that $a+b=c$ and $\operatorname{rad}(abc) < c^{1-\epsilon}$ for any $\epsilon > 0$. Using the optimized methods in a recent work of Browning, Lichtman and Teräväinen, we showed that the number of those triples with $c \leqslant X$ is $O\left(X^{56/85+\varepsilon}\right)$ for any $\varepsilon > 0$, where $\frac{56}{85} \approx 0.658824$. This constitutes an improvement of the previous bound $O\left(X^{33/50}\right)$.
\end{abstract}

\maketitle

\tableofcontents

\section{Introduction}
Let $n$ denotes a positive integer, $p$ denotes a prime and write
\begin{equation}
\operatorname{rad}(n) = \prod_{p | n} p.
\end{equation}
We say a triple $(a,b,c)$ of coprime positive integers $a, b, c$ is an $abc$ triple of exponent $\lambda$ if
$$
a+b=c \quad \text{and} \quad \operatorname{rad}(abc) < c^{\lambda}.
$$
The famous $abc$ conjecture, proposed by Masser and Oesterlé, asserts that there are only finitely many $abc$ triples of exponent $\lambda$ for any $\lambda < 1$. Now the best result in this direction is due to Stewart and Yu, who showed that there are finitely many $abc$ triples satisfy
\begin{equation}
\operatorname{rad}(abc) < (\log c)^{3 - \epsilon}.
\end{equation}
For more historical progress of the $abc$ conjecture, we refer the readers to \cite{abc66}.

Now, we are focusing on the exceptional set in the $abc$ conjecture. We first define $N_{\lambda}(X)$ as the number of $abc$ triples of exponent $\lambda$ in $[1, X]^3$ as $X \to \infty$. A "trivial" bound states that
\begin{theorem}\label{t1} ("Trivial" bound).
Let $\lambda > 0$. Then we have
$$
N_{\lambda}(X) \ll x^{\frac{2}{3}\lambda + \epsilon}
$$
for any $\epsilon > 0$.
\end{theorem}
For the proof, one can see Lichtman's recent note \cite{abcNoteLichtman}.

In 2024, Browning, Lichtman and Teräväinen \cite{abc66} developed a system of combinatorial bounds and improved Theorem~\ref{t1}. Their result is the first power--saving improvement over the "trivial" bound for $\lambda$ close to 1 (actually, for $0.99 < \lambda < 1.001$).
\begin{theorem}\label{t2} (Browning--Lichtman--Teräväinen bound).
Let $0 < \lambda < 1.001$. Then we have
$$
N_{\lambda}(X) \ll x^{\frac{33}{50}} = x^{0.66}.
$$
\end{theorem}

In the present paper, we use the method of Browning, Lichtman and Teräväinen \cite{abc66} to improve their result and show that without further optimization, the best exponent their current method can reach is $\frac{56}{85}$.
\begin{theorem}\label{t3}
For any $\varepsilon > 0$, there exists a positive constant $\delta = \delta(\varepsilon)$ such that for $0 < \lambda < 1 + \delta$, we have
$$
N_{\lambda}(X) \ll x^{\frac{56}{85} + \varepsilon}.
$$
\end{theorem}

In this paper, we put $\varepsilon > 0$, $0 < \delta < 10^{-100}$ and $\theta = \frac{56}{85} + \varepsilon$. We also suppose that $\epsilon$ is a sufficiently small positive number.

\section{Number of solutions to Diophantine equations}
We define a counting function $S_{\alpha, \beta, \gamma}(X)$ for $\alpha, \beta, \gamma > 0$ as the same as in \cite{abc66}: $S_{\alpha, \beta, \gamma}(X)$ denotes the number of $(a, b, c) \in \mathbb{N}^3$ with $\operatorname{gcd}(a, b, c) = 1$ such that
$$
a, b, c \in [1, X], \quad a+b=c, \quad \operatorname{rad}(a) \leqslant a^{\alpha}, \quad \operatorname{rad}(b) \leqslant b^{\beta}, \quad \operatorname{rad}(c) \leqslant c^{\gamma}.
$$
Then we have
\begin{equation}
N_{\lambda}(X) \leqslant \max_{\substack{\alpha, \beta, \gamma > 0 \\ \alpha + \beta + \gamma \leqslant \lambda}} S_{\alpha, \beta, \gamma}(X).
\end{equation}
We shall use a standard dyadic decomposition to define a variant of $S_{\alpha, \beta, \gamma}(X)$: Let $S^{*}_{\alpha, \beta, \gamma}(X)$ denotes the number of $(a, b, c) \in \mathbb{N}^3$ with $\operatorname{gcd}(a, b, c) = 1$ such that
$$
c \in \left[\frac{X}{2}, X \right], \quad a+b=c, \quad \operatorname{rad}(a) \sim a^{\alpha}, \quad \operatorname{rad}(b) \sim b^{\beta}, \quad \operatorname{rad}(c) \sim c^{\gamma}.
$$
Then, by the pigeonhole principle we have
\begin{equation}
S_{\alpha, \beta, \gamma}(X) \ll (\log X)^4 \max_{\substack{\alpha^{\prime} \leqslant \alpha \\ \beta^{\prime} \leqslant \beta \\ \gamma^{\prime} \leqslant \gamma}} \max_{Y \in [1, X]} S^{*}_{\alpha^{\prime}, \beta^{\prime}, \gamma^{\prime}}(Y).
\end{equation}
Hence, in order to prove Theorem~\ref{t3}, we only need to show that
\begin{equation}
S^{*}_{\alpha, \beta, \gamma}(X) \ll X^{\theta} (\log X)^{-4}.
\end{equation}
We need the following important lemma to reduce the problem into bounding the number of solutions to some Diophantine equations.

\begin{lemma}\label{l21} ([\cite{abc66}, Proposition 2.1]).
Let $\alpha, \beta, \gamma \in (0, 1]$ be fixed and let $X \geqslant 2$. For any $\epsilon > 0$ there exists an integer $d = d(\epsilon) \geqslant 1$ such that the following holds: There exist $X_1, \ldots, X_d, Y_1, \ldots, Y_d, Z_1, \ldots, Z_d \geqslant 1$ satisfying
$$
X^{\alpha - \epsilon} \leqslant \prod_{1 \leqslant j \leqslant d}X_j \leqslant X^{\alpha + \epsilon}, \quad X^{\beta - \epsilon} \leqslant \prod_{1 \leqslant j \leqslant d}Y_j \leqslant X^{\beta + \epsilon}, \quad X^{\gamma - \epsilon} \leqslant \prod_{1 \leqslant j \leqslant d}Z_j \leqslant X^{\gamma + \epsilon},
$$
$$
\prod_{1 \leqslant j \leqslant d}X_j^j \leqslant X, \quad \prod_{1 \leqslant j \leqslant d}Y_j^j \leqslant X, \quad X^{1-\epsilon^2} \leqslant \prod_{1 \leqslant j \leqslant d}Z_j^j \leqslant X,
$$
and pairwise coprime integers $1 \leqslant c_1, c_2, c_3 \leqslant X^{\epsilon}$, such that
$$
S^{*}_{\alpha, \beta, \gamma}(X) \ll X^{\epsilon} B_d(\mathbf{c}, \mathbf{X}, \mathbf{Y}, \mathbf{Z}),
$$
where
\begin{align}
\nonumber B_d(\mathbf{c}, \mathbf{X}, \mathbf{Y}, \mathbf{Z}) = \# & \left\{ (\mathbf{x}, \mathbf{y}, \mathbf{z}) \in \mathbb{N}^{3d} : x_i \sim X_i,\ y_i \sim Y_i,\ z_i \sim Z_i, \right. \\
\nonumber & c_1 \prod_{j \leqslant d}x_j^j + c_2 \prod_{j \leqslant d}y_j^j = c_3 \prod_{j \leqslant d}z_j^j, \\
\nonumber & \left. \operatorname{gcd}\left(c_1 \prod_{j \leqslant d}x_j,\ c_2 \prod_{j \leqslant d}y_j,\ c_3 \prod_{j \leqslant d}z_j \right) = 1 \right\}
\end{align}
for $\mathbf{c} \in \mathbb{Z}^3$ and $\mathbf{X}, \mathbf{Y}, \mathbf{Z} \in \mathbb{R}^{d}_{>0}$.
\end{lemma}

Now we give some upper bounds for the integer points $B_d(\mathbf{c}, \mathbf{X}, \mathbf{Y}, \mathbf{Z})$. These lemmas are proved in \cite{abc66} and they will be used to give combinatorial bounds for $\nu$ in next section.

\begin{lemma}\label{l22} (Fourier bound, [\cite{abc66}, Proposition 3.1]).
Let $d \geqslant 1$, $\epsilon > 0$ and $A \geqslant 1$ be fixed. Let
$$
X_1, \ldots, X_d, Y_1, \ldots, Y_d, Z_1, \ldots, Z_d \geqslant 1.
$$
Let $\mathbf{c} = (c_1, c_2, c_3) \in \mathbb{Z}^3$ satisfy $0 < |c_1|, |c_2|, |c_3| \leqslant \max_{1 \leqslant i \leqslant d}(X_i Y_i Z_i)^A$. Then we have
$$
B_d(\mathbf{c}, \mathbf{X}, \mathbf{Y}, \mathbf{Z}) \ll \max_{1 \leqslant i \leqslant d}(X_i Y_i Z_i)^{\epsilon} \frac{\prod_{j \leqslant d} \left(X_j Y_j Z_j (Y_j + Z_j)\right)^{\frac{1}{2}} }{\max_{i = 1} \prod_{j \equiv 0 (\bmod i)} Z_j^{\frac{1}{2}}}.
$$
\end{lemma}

\begin{lemma}\label{l23} (Geometry bound, [\cite{abc66}, Proposition 3.2]).
Let $d \geqslant 1$ and $\epsilon > 0$ be fixed. Let
$$
X_1, \ldots, X_d, Y_1, \ldots, Y_d, Z_1, \ldots, Z_d \geqslant 1.
$$
Let $\mathbf{c} = (c_1, c_2, c_3) \in \mathbb{Z}^3$ have non--zero and pairwise coprime coordinates. Then we have
$$
B_d(\mathbf{c}, \mathbf{X}, \mathbf{Y}, \mathbf{Z}) \ll \max_{1 \leqslant i \leqslant d}(X_i Y_i Z_i)^{\epsilon} \min_{I, I^{\prime}, I^{\prime \prime} \subset [d]} \left(\prod_{i \in I}X_i \prod_{i \in I^{\prime}}Y_i \prod_{i \in I^{\prime \prime}}Z_i \right) \left(1+\frac{\prod_{i \notin I}X_i^i \prod_{i \notin I^{\prime}}Y_i^i \prod_{i \notin I^{\prime \prime}}Z_i^i}{\max\left(|c_1|\prod_{i}X_i^i, |c_2|\prod_{i}Y_i^i, |c_3|\prod_{i}Z_i^i\right)} \right).
$$
\end{lemma}

\begin{lemma}\label{l24} (Determinant bound, [\cite{abc66}, Proposition 3.5]).
Let $d \geqslant 1$ and let
$$
X_1, \ldots, X_d, Y_1, \ldots, Y_d, Z_1, \ldots, Z_d \geqslant 1.
$$
Let $\mathbf{c} = (c_1, c_2, c_3) \in \mathbb{Z}_{\neq 0}^3$. Then we have
$$
B_d(\mathbf{c}, \mathbf{X}, \mathbf{Y}, \mathbf{Z}) \ll \max_{1 \leqslant i \leqslant d}(X_i Y_i Z_i)^{\epsilon} \prod_{i \leqslant d}(X_i Y_i Z_i) \min_{p,q \geqslant 1}\left((X_p Y_q)^{-1} \min\left(X_p^{1/q} Y_q^{1/p}\right) \right).
$$
\end{lemma}

\begin{lemma}\label{l25} (Thue bound, [\cite{abc66}, Proposition 3.6]).
Let $d \geqslant 1$ and let
$$
X_1, \ldots, X_d, Y_1, \ldots, Y_d, Z_1, \ldots, Z_d \geqslant 1.
$$
Let $\mathbf{c} = (c_1, c_2, c_3) \in \mathbb{Z}_{\neq 0}^3$. Then we have
$$
B_d(\mathbf{c}, \mathbf{X}, \mathbf{Y}, \mathbf{Z}) \ll \max_{1 \leqslant i \leqslant d}(X_i Y_i Z_i)^{\epsilon} \prod_{i \leqslant d}(X_i Y_i Z_i) \min_{p \geqslant 2}\left(\prod_{\substack{j \leqslant d \\ p | j }}(X_j Y_j)^{-1} \right).
$$
\end{lemma}

\section{Upper bounds for $\nu$}
In this section we shall use all things proved above to bound $B_d(\mathbf{c}, \mathbf{X}, \mathbf{Y}, \mathbf{Z})$ for any pairwise coprime integers $1 \leqslant |c_1|, |c_2|, |c_3| \leqslant X^{\epsilon^2}$, any fixed $d \geqslant 1$ and any choice of $X_i, Y_i, Z_i \geqslant 1$ for $1 \leqslant i \leqslant d$ that satisfies conditions in Lemma~\ref{l21}. Moreover, we have 
\begin{equation}
\alpha + \beta + \gamma \leqslant \lambda \leqslant 1 + \delta - \epsilon.
\end{equation}

We define $a_i, b_i, c_i$ by writing
\begin{equation}
X_i = X^{a_i}, \quad Y_i = X^{b_i}, \quad Z_i = X^{c_i}
\end{equation}
for $i \leqslant d$ and $a_i = b_i = c_i = 0$ for $i > d$. We write $s_i = a_i + b_i + c_i$. By the conditions in Lemma~\ref{l21}, we can assume that
\begin{equation}
\sum_{i \leqslant d} i a_i,\ \sum_{i \leqslant d} i b_i \leqslant 1, \quad 1 - \epsilon^2 \leqslant \sum_{i \leqslant d} i c_i \leqslant 1.
\end{equation}
By (6) and [\cite{abc66}, (1.2)], we can also assume that
\begin{equation}
\sum_{i \leqslant d}(a_i + b_i),\ \sum_{i \leqslant d}(a_i + c_i),\ \sum_{i \leqslant d}(b_i + c_i) \geqslant \theta - \epsilon^2
\end{equation}
and
\begin{equation}
\sum_{i \leqslant d} s_i \leqslant 1 + \delta - \epsilon.
\end{equation}

We define
\begin{equation}
\nu = \frac{\log B_d(\mathbf{c}, \mathbf{X}, \mathbf{Y}, \mathbf{Z})}{\log X} + 2 \epsilon^2.
\end{equation}
Then we only need to show that
\begin{equation}
\nu \leqslant \theta.
\end{equation}

Now we shall rewrite Lemmas~\ref{l22}--\ref{l25} in terms of an upper bound for $\nu$ using parameters $a_i$, $b_i$, $c_i$.
\begin{lemma}\label{Fourier} (Fourier bound).
We have
$$
\nu \leqslant \frac{1}{2}\left(1+\delta + \sum_{i \leqslant d} \max(a_i, b_i) - \max_{m > 1}(a_m, b_m) \right).
$$
\end{lemma}

\begin{lemma}\label{Geometry} (Geometry bound).
We have
$$
\nu \leq \delta + \min_{I, I^{\prime}, I^{\prime \prime} \subset [d]} \left( \max\left(1, \sum_{i \in I}i a_i + \sum_{i \in I^{\prime}}i b_i + \sum_{i \in I^{\prime \prime}}i c_i\right) - \sum_{i \in I} a_i - \sum_{i \in I^{\prime}} b_i - \sum_{i \in I^{\prime \prime}} c_i \right)
$$
or
$$
\nu \leq 4 \epsilon^2 + \min_{I, I^{\prime}, I^{\prime \prime} \subset [d]} \left( \sum_{i \notin I} a_i + \sum_{i \notin I^{\prime}} b_i + \sum_{i \notin I^{\prime \prime}} c_i + \max\left(0, \sum_{i \in I}i a_i + \sum_{i \in I^{\prime}}i b_i + \sum_{i \in I^{\prime \prime}}i c_i - 1\right) \right).
$$
\end{lemma}

\begin{lemma}\label{Determinant} (Determinant bound).
We have
$$
\nu < \min_{p, q \geqslant 1}\left(1+\delta - a_p - b_q + \min\left(\frac{a_p}{q}, \frac{b_q}{p} \right) \right).
$$
\end{lemma}

\begin{lemma}\label{Thue} (Thue bound).
We have
$$
\nu < 1+\delta - \max_{p \geqslant 2}\left(\sum_{p | i}(a_i + b_i) \right).
$$
\end{lemma}

We first show that we can assume that
\begin{equation}
2 \theta - 1 - \delta \leqslant \sum_{i \leqslant d} a_i,\ \sum_{i \leqslant d} b_i,\ \sum_{i \leqslant d} c_i \leqslant 1 - \theta + \delta - \frac{1}{2}\epsilon.
\end{equation}
If $\sum_{i \leqslant d} c_i > 1 - \theta + \delta - \frac{1}{2}\epsilon$, then we have
\begin{equation}
\sum_{i \leqslant d} (a_i + b_i) \leqslant \theta - \frac{1}{2}\epsilon.
\end{equation}
By [\cite{abc66}, (1.2)], Theorem~\ref{t3} is proved. If $\sum_{i \leqslant d} c_i < 2 \theta - 1 - \delta$ and all of the three sums are $\leqslant 1 - \theta + \delta - \frac{1}{2}\epsilon$, then we have
\begin{equation}
\sum_{i \leqslant d} (b_i + c_i) \leqslant (2 \theta - 1 - \delta) + \left(1 - \theta + \delta - \frac{1}{2}\epsilon \right) = \theta - \frac{1}{2}\epsilon.
\end{equation}
Again, Theorem~\ref{t3} is proved by [\cite{abc66}, (1.2)].

\section{Proof of Theorem 1.1}
From now on, we ignore the presence of $\delta$ and $\epsilon$ in many places, since all the contributions of them can be bounded by $\varepsilon$. We define the parameters $\delta_a, \delta_b, \delta_c, \delta_{ab}, \delta_{ac}, \delta_{bc}, \delta_s$ by
\begin{equation}
\delta_a = \frac{1}{3} - \sum_{i \leqslant d} a_i, \quad \delta_b = \frac{1}{3} - \sum_{i \leqslant d} b_i, \quad \delta_c = \frac{1}{3} - \sum_{i \leqslant d} c_i,
\end{equation}
\begin{equation}
\delta_{ab} = \delta_a + \delta_b, \quad \delta_{ac} = \delta_a + \delta_c, \quad \delta_{bc} = \delta_b + \delta_c, \quad \delta_s = \delta_a + \delta_b + \delta_c.
\end{equation}
By (9), (16) and (17) we know that
\begin{equation}
\delta_{ab}, \delta_{ac}, \delta_{bc} \leqslant \frac{2}{3} - \theta.
\end{equation}
By (13) and (16) we have
\begin{equation}
\theta - \frac{2}{3} \leqslant \delta_a, \delta_b, \delta_c \leqslant \frac{4}{3} - 2 \theta.
\end{equation}
By (10) and (16) we know that
\begin{equation}
1 - \delta_s \leqslant 1 + \delta
\end{equation}
and
\begin{equation}
2 \delta_s = \delta_{ab} + \delta_{ac} + \delta_{bc} \leqslant 2 - 3 \theta.
\end{equation}
Then, by (20) and (21) we have
\begin{equation}
- \delta < \delta_s \leqslant 1 - \frac{3}{2} \theta.
\end{equation}
Note that these inequalities
\begin{equation}
\sum_{i \geqslant 2}(i-1)a_i \leqslant \frac{2}{3} + \delta_a, \quad \sum_{i \geqslant 3}(i-2)a_i \leqslant \frac{1}{3} + a_1 + 2 \delta_a, \quad \sum_{i \geqslant 4}(i-3)a_i \leqslant 2 a_1 + a_2 + 3 \delta_a 
\end{equation}
follow by (8) and subtracting. Similar inequalities hold for $b_i$ and $c_i$. By Lemma~\ref{Thue}, we know that
$$
\nu \leqslant 1 + \delta - \max_{p \geqslant 2} \sum_{p \mid i}(a_i + b_i)
$$
and similar results hold for $a_i + c_i$ and $b_i + c_i$. Thus we can assume that
\begin{equation}
a_i + b_i, a_i + c_i, b_i + c_i < 1 - \theta
\end{equation}
for every $i \geqslant 2$. Moreover, we can assume that
\begin{equation}
a_2 + b_2 + a_4 + b_4, a_2 + c_2 + a_4 + c_4, b_2 + c_2 + b_4 + c_4, < 1 - \theta.
\end{equation}
Now, (24) and (25) imply that
\begin{equation}
s_2 + s_4, s_3, s_5 \leqslant \frac{3}{2}(1 - \theta) = \frac{3}{2} - \frac{3}{2} \theta.
\end{equation}
By (16) and (23) we also know that
\begin{equation}
\sum_{i \geqslant 1} s_i = 1 - \delta_s, \quad \sum_{i \geqslant 2}(i-1) s_i \leqslant 2 + \delta_s, \quad \sum_{i \geqslant 3}(i-2) s_i \leqslant 1 + s_1 + 2 \delta_s, \quad \sum_{i \geqslant 4}(i-3) s_i \leqslant 2 s_1 + s_2 + 3 \delta_s.
\end{equation}

If $s_1 + s_2 > 1 - \theta$, then by Lemma~\ref{Geometry} and (26), we have
\begin{align}
\nonumber \nu \leqslant&\ \max(1, s_1 + 2 s_2) - s_1 - s_2 + \delta \\
\nonumber =&\ \max(1 - s_1 - s_2, s_2) + \delta \\
<&\ \max\left(\theta, \frac{3}{2} - \frac{3}{2} \theta\right) = \theta
\end{align}
since $\theta > 0.6$. Now we can assume that $s_1 + s_2 \leqslant 1 - \theta$.

For any $i \geqslant 3$, let $\tau_i$ be an element in $\{a_i, b_i, c_i, a_i + b_i, a_i + c_i, b_i + c_i, s_i \}$. By Lemma~\ref{Geometry} we know that
\begin{align}
\nonumber \nu \leqslant&\ \max(1, s_1 + 2 s_2 + i \tau_i) - s_1 - s_2 - \tau_i + \delta \\
=&\ \max(1 - s_1 - s_2 - \tau_i, s_2 + (i-1) \tau_i ) + \delta 
\end{align}
and
\begin{align}
\nonumber \nu \leqslant&\ \max(1, s_1 + 3 \tau_3) - s_1 - \tau_3 + \delta \\
=&\ \max(1 - s_1 - \tau_3, 2 \tau_3 ) + \delta. 
\end{align}
Combining (29) and (30), we know that $\nu \leqslant \theta$ if
\begin{equation}
\tau_3 \in \left(1-\theta -s_1 -s_2,\ \frac{1}{2}\theta - \frac{1}{2}s_2 \right) \cup \left(1-\theta -s_1,\ \frac{1}{2}\theta  \right).
\end{equation}

By (23) we know that
\begin{equation}
\sum_{i \geqslant 4}a_i \leqslant \sum_{i \geqslant 4}(i-3) a_i \leqslant 2 a_1 + a_2 + 3 \delta_a
\end{equation}
and
\begin{equation}
\sum_{i \geqslant 5}a_i \leqslant \frac{1}{2}\left( \sum_{i \geqslant 4}(i-3) a_i - a_4 \right) \leqslant \frac{1}{2}\left( 2 a_1 + a_2 - a_4 + 3 \delta_a \right).
\end{equation}
By (16), these imply that
\begin{equation}
a_3 = \frac{1}{3} - \delta_a - a_1 - a_2 - \sum_{i \geqslant 4}a_i \geqslant \frac{1}{3} - 3 a_1 - 2 a_2 - 4 \delta_a
\end{equation}
and
\begin{equation}
a_3 = \frac{1}{3} - \delta_a - a_1 - a_2 - a_4 - \sum_{i \geqslant 5}a_i \geqslant \frac{1}{3} - 2 a_1 - \frac{3}{2} a_2 - \frac{1}{2} a_4 - \frac{5}{2} \delta_a.
\end{equation}
Note that (34) and (35) also hold for $b_3$ and $c_3$. Adding up these corresponding lower bounds, we have
\begin{equation}
s_3 \geqslant 1 - 3 s_1 - 2 s_2 - 4 \delta_s
\end{equation}
and
\begin{equation}
s_3 \geqslant 1 - 2 s_1 - \frac{3}{2} s_2 - \frac{1}{2} s_4 - \frac{5}{2} \delta_s.
\end{equation}
Now, we split the argument according to whether $s_2 \geqslant k$ or $s_2 < k$, where
\begin{equation}
k = \frac{49}{12} - \frac{23}{4}\theta \approx 0.2951.
\end{equation}
Without loss of generality, we shall assume that $a_3 \geqslant b_3 \geqslant c_3$ in all that follows.

\subsection{Case 1: $s_2 \geqslant k$}
By the assumption $s_1 + s_2 \leqslant 1 - \theta$ we know that
\begin{equation}
s_1 \leqslant 1 - \theta - s_2 \leqslant 1 - \theta - k.
\end{equation}
By (26) we know that
\begin{equation}
s_4 \leqslant \frac{3}{2} - \frac{3}{2}\theta - s_2 \leqslant \frac{3}{2} - \frac{3}{2}\theta - k.
\end{equation}

\subsubsection{Subcase 1.1: $b_3 \leqslant 1 - \theta - s_1 - s_2$} Because $c_3 \leqslant b_3$, we have
\begin{align}
\nonumber b_3 + c_3 \leqslant 2 b_3 \leqslant&\ 2 (1 - \theta - s_1 - s_2) \\
\nonumber =&\ 2 - 2 \theta - 2 s_1 - 2 s_2 \\
\leqslant&\ 2 - 2 \theta - 2 s_2.
\end{align}
Note that we have
\begin{equation}
2 - 2 \theta - 2 s_2 \leqslant \frac{1}{2}\theta - \frac{1}{2}s_2
\end{equation}
since $s_2 \geqslant \frac{4}{3} - \frac{5}{3}\theta$. We also have
\begin{equation}
\frac{4}{3} - \frac{5}{3}\theta \leqslant k = \frac{49}{12} - \frac{23}{4}\theta.
\end{equation}
Since $s_2 \geqslant k$ in this case, we have
\begin{equation}
b_3 + c_3 \leqslant \frac{1}{2}\theta - \frac{1}{2}s_2.
\end{equation}
If $b_3 + c_3$ is in the interval (31), we get $\nu \leqslant \theta$. Otherwise we must have
\begin{equation}
b_3 + c_3 \leqslant 1-\theta -s_1 -s_2.
\end{equation}
(We will repeat similar discussions for many times in the following.) Now, by (35), (45) and (39) we can lower bound $a_3$ by
\begin{align}
\nonumber a_3 = s_3 - (b_3 + c_3) \geqslant&\ \left( 1 - 2 s_1 - \frac{3}{2} s_2 - \frac{1}{2} s_4 - \frac{5}{2} \delta_a \right) - (1-\theta -s_1 -s_2) \\
\nonumber =&\ \theta - \frac{5}{2} \delta_a - s_1 - \frac{1}{2} \left(s_2 + s_4\right) \\
\nonumber \geqslant&\ \theta - \frac{5}{2} \delta_a - (1 - \theta - k) - \frac{1}{2} \left(\frac{3}{2} - \frac{3}{2} \theta \right) \\
=&\ \frac{11}{4} \theta - \frac{5}{2} \delta_a - \frac{7}{4} + k.
\end{align}
Now, (46) and (16) ensure that
\begin{equation}
\frac{11}{4} \theta - \frac{5}{2} \delta_a - \frac{7}{4} + k \leqslant a_3 \leqslant \frac{1}{3} - \delta_a.
\end{equation}
By (47) and (19), we know that
\begin{align}
\nonumber \frac{11}{4} \theta + k - \frac{7}{4} - \frac{1}{3} \leqslant&\ \frac{3}{2} \delta_a \\
\nonumber \leqslant&\ \frac{3}{2} \left(\frac{4}{3} -2 \theta \right) \\
\nonumber \frac{11}{4} \theta + k - \frac{25}{12} \leqslant&\ 2 - 3 \theta \\
k \leqslant&\ \frac{49}{12} - \frac{23}{4}\theta.
\end{align}
Now (48) contradicts our assumption since the contributions of $\delta$ are omitted.

\subsubsection{Subcase 1.2: $b_3 > 1 - \theta - s_1 - s_2$} By (31) and a similar discussion as in (44)--(45), we can assume that
\begin{equation}
b_3 \geqslant \frac{1}{2}\theta - \frac{1}{2}s_2.
\end{equation}
By (23), we have
\begin{equation}
\sum_{i \geqslant 4}(i-2)b_i = \sum_{i \geqslant 3}(i-2)b_i - b_3 \leqslant \frac{1}{3} + b_1 - b_3 + 2 \delta_b.
\end{equation}
We also have
\begin{equation}
b_1 \leqslant s_1 \leqslant 1 - \theta - s_2.
\end{equation}
Thus, by (50), (51) and (49) we have
\begin{align}
\nonumber \sum_{i \geqslant 4}b_i \leqslant&\ \frac{1}{2} \sum_{i \geqslant 4}(i-2)b_i \\
\nonumber \leqslant&\ \frac{1}{2} \left(\frac{1}{3} + b_1 - b_3 + 2 \delta_b \right) \\
\nonumber \leqslant&\ \frac{1}{2} \left(\frac{1}{3} + (1 - \theta - s_2) - \left(\frac{1}{2}\theta - \frac{1}{2}s_2 \right) + 2 \left(\frac{4}{3} -2 \theta \right) \right) \\
\nonumber =&\ 2 - \frac{11}{4}\theta - \frac{1}{4}s_2.
\end{align}
We want to show that
\begin{equation}
2 - \frac{11}{4}\theta - \frac{1}{4}s_2 < \frac{1}{2}\theta - \frac{1}{2}s_2.
\end{equation}
Note that this is equivalent to
\begin{equation}
s_2 < 13 \theta - 8.
\end{equation}
Now, by the assumption above we know that $s_2 < 1 - \theta$, and we have
\begin{equation}
1 - \theta < 13 \theta - 8
\end{equation}
since $\theta > \frac{9}{14} \approx 0.6428$. Combining (52)--(54) we know that
\begin{equation}
\sum_{i \geqslant 4}a_i,\ \sum_{i \geqslant 4}b_i \leqslant \frac{1}{2}\theta - \frac{1}{2}s_2
\end{equation}
and by (31) and a similar discussion as in (44)--(45) we can assume
\begin{equation}
a_4, \ b_4, \ a_5, \ b_5, \ a_6, \ b_6 \leqslant 1-\theta -s_1 -s_2.
\end{equation}

Now, by Lemma~\ref{Fourier} we have
\begin{equation}
\nu < \frac{1}{2}\left(1+\delta +\sum_{i \leqslant d}\max(a_i, b_i) - \max(a_2, b_2) \right).
\end{equation}
Using (16), this implies that
\begin{align}
\nonumber 2 \nu -1-\delta < \sum_{i \neq 2}\max(a_i, b_i) \leqslant&\ \sum_{i \neq 2, \ i \leqslant 6}\max(a_i, b_i) + \sum_{i \geqslant 7}(a_i + b_i) \\
\nonumber =&\ \sum_{i \neq 2, \ i \leqslant 6}\max(a_i, b_i) + \frac{2}{3} - \delta_{ab} - \sum_{i \leqslant 6}(a_i + b_i) \\
=&\ \frac{2}{3} - \delta_{ab} - \sum_{i \neq 2, \ i \leqslant 6}\min(a_i, b_i) - (a_2 + b_2).
\end{align}

We then give a lower bound for $a_2 + b_2$. By (23), we have
\begin{equation}
4 \sum_{i \geqslant 7}a_i \leqslant \sum_{i \geqslant 7}(i - 3)a_i = \left( \sum_{i \geqslant 4}(i - 3)a_i \right) - a_4 - 2 a_5 - 3 a_6 = (2 a_1 + a_2 + 3 \delta_a) - a_4 - 2 a_5 - 3 a_6,
\end{equation}
whence
\begin{equation}
\frac{1}{3} - \delta_a = \sum_{i \leqslant 6}a_i + \sum_{i \geqslant 7}a_i \leqslant \sum_{i \leqslant 6}a_i + \frac{1}{4} (2 a_1 + a_2 + 3 \delta_a - a_4 - 2 a_5 - 3 a_6) = \frac{1}{4} \sum_{i \leqslant 6}(7 - i) a_i + \frac{3}{4} \delta_a.
\end{equation}
Then we have
\begin{equation}
a_2 \geqslant \frac{4}{15} - \frac{1}{5} \sum_{i \neq 2, \ i \leqslant 6}(7 - i) a_i - \frac{7}{5} \delta_a
\end{equation}
and
\begin{equation}
b_2 \geqslant \frac{4}{15} - \frac{1}{5} \sum_{i \neq 2, \ i \leqslant 6}(7 - i) b_i - \frac{7}{5} \delta_b.
\end{equation}
Since $\min(a_3, b_3) = b_3$, we now have
\begin{align}
\nonumber 2 \nu -1-\delta <&\ \frac{2}{3} - \delta_{ab} - \sum_{i \neq 2, \ i \leqslant 6}\min(a_i, b_i) - \left( \frac{8}{15} - \frac{1}{5} \sum_{i \neq 2, \ i \leqslant 6}(7 - i) (a_i + b_i) - \frac{7}{5} \delta_{ab} \right) \\
\nonumber \leqslant&\ \frac{2}{15} + \frac{2}{5} \delta_{ab} + \frac{1}{5} \left( 6 \max(a_1, b_1) + \min(a_1, b_1) + 4 a_3 - b_3 \right. \\
& \left. \qquad \qquad \qquad \quad + 3 \max(a_4, b_4) + 2 \max(a_5, b_5) + \max(a_6, b_6) \right).
\end{align}
Using (24) and (49), we have
\begin{align}
\nonumber 4 a_3 - b_3 \leqslant&\ 4 (1 - \theta - b_3) - b_3 \\
\nonumber <&\ 4 \left(1 - \theta - \left(\frac{1}{2}\theta - \frac{1}{2}s_2 \right) \right) - \left(\frac{1}{2}\theta - \frac{1}{2}s_2 \right) \\
=&\ \frac{5}{2} s_2 + 4 - \frac{13}{2} \theta.
\end{align}
Finally, by (63)--(64) we have
\begin{align}
\nonumber 2 \nu -1-\delta <&\ \frac{2}{15} + \frac{2}{5} \delta_{ab} + \frac{1}{5} \left( 6 \max(a_1, b_1) + \min(a_1, b_1) + 4 a_3 - b_3 \right. \\
\nonumber & \left. \qquad \qquad \qquad \quad + 3 \max(a_4, b_4) + 2 \max(a_5, b_5) + \max(a_6, b_6) \right) \\
\nonumber <&\ \frac{2}{15} + \frac{2}{5} \delta_{ab} + \frac{1}{5} \left( 6 s_1 + \left( \frac{5}{2} s_2 + 4 - \frac{13}{2} \theta \right) + 6 (1 - \theta - s_1 - s_2) \right) \\
\nonumber <&\ \frac{2}{15} + \frac{2}{5} \delta_{ab} + \frac{1}{5} \left( 6 s_1 + \frac{5}{2} s_2 + 4 - \frac{13}{2} \theta  + 6 - 6 \theta - 6 s_1 - 6 s_2 \right) \\
\nonumber =&\ \frac{32}{15} - \frac{5}{2} \theta - \frac{7}{10}s_2 + \frac{2}{5} \delta_{ab} \\
\leqslant&\ \frac{32}{15} - \frac{5}{2} \theta - \frac{7}{10} \left(\frac{49}{12} - \frac{23}{4}\theta \right) + \frac{2}{5} \left(\frac{2}{3} -\theta \right) < 0.283, \\
\nu <&\ \frac{1}{2}(1+0.283) < 0.65.
\end{align}

\subsection{Case 2: $s_2 < k$}
By (24) we have
\begin{equation}
2 b_3 \leqslant a_3 + b_3 < 1 - \theta,
\end{equation}
so that
\begin{equation}
b_3 \leqslant \frac{1}{2} - \frac{1}{2}\theta.
\end{equation}
We want to show that
\begin{equation}
\frac{1}{2} - \frac{1}{2}\theta \leqslant \frac{1}{2}\theta - \frac{1}{2}s_2.
\end{equation}
Note that (69) holds if
\begin{equation}
s_2 \leqslant 2 \theta -1.
\end{equation}
Because we have
\begin{equation}
s_2 < k = \frac{49}{12} - \frac{23}{4}\theta \leqslant 2 \theta - 1
\end{equation}
when $\theta \geqslant \frac{61}{93} \approx 0.6559$, we deduce that
\begin{equation}
c_3 \leqslant b_3 < \frac{1}{2}\theta - \frac{1}{2}s_2.
\end{equation}
By (31) and a similar discussion as in (44)--(45), we can assume that
\begin{equation}
c_3 \leqslant b_3 \leqslant 1 - \theta - s_1 - s_2.
\end{equation}
Then (36) gives that
\begin{align}
\nonumber a_3 = s_3 - (b_3 + c_3) \geqslant&\ (1 - 3 s_1 - 2 s_2 - 4 \delta_s) - 2 (1 - \theta - s_1 - s_2) \\
=&\ 2 \theta - 1 - s_1 - 4 \delta_s.
\end{align}

We first prove the bound (12) in two cases:
$$
a_3 \geqslant 2 \theta - 1 \qquad \text{and} \qquad b_3 + c_3 < \frac{1}{2} \theta - \frac{1}{2} s_2.
$$

\subsubsection{Subcase 2.1: $a_3 \geqslant 2 \theta - 1$} In this case we have
\begin{equation}
b_3, \ c_3 \leqslant 1 - \theta - a_3 \leqslant 2 - 3 \theta.
\end{equation}
Let
\begin{equation}
M = \max_{i \geqslant 4} \max(b_i, c_i).
\end{equation}
If $M > 3 - \frac{9}{2} \theta$, by Lemma~\ref{Determinant} we know that
\begin{equation}
\nu = 1 + \delta - a_3 - M + \min\left(\frac{M}{3}, \frac{a_3}{4} \right) \leqslant 1 + \delta - a_3 - \frac{2}{3} M \leqslant 1 - (2 \theta - 1) - (2 - 3 \theta) = \theta.
\end{equation}
Thus we can assume that $\max(b_i, c_i) \leqslant 3 - \frac{9}{2} \theta$ for $i \geqslant 4$. Then we have
\begin{equation}
b_i + c_i \leqslant 2 \max(b_i, c_i) \leqslant 6 - 9 \theta
\end{equation}
for $i \geqslant 4$. Moreover, by (8) and (16) we know that
\begin{equation}
\sum_{i \leqslant d}(i - 1) (b_i + c_i) \leqslant \frac{4}{3} + \delta_{bc}.
\end{equation}
Using the second form of Lemma~\ref{Geometry}, we have
\begin{align}
\nonumber \nu \leqslant&\ \epsilon + a_3 + b_3 + \min(b_4, c_4) + \sum_{i \geqslant 5}(b_i + c_i) \\
& + \max\left(0,\ \sum_{i} i s_i - 3(a_3 + b_3) - 4 \min(b_4, c_4) - \sum_{i \geqslant 5} i (b_i + c_i) - 1 \right).
\end{align}
Now, define
\begin{equation}
\nu_1 = a_3 + b_3 + \min(b_4, c_4) + \sum_{i \geqslant 5}(b_i + c_i)
\end{equation}
and
\begin{equation}
\nu_2 = \sum_{i} i s_i - 2(a_3 + b_3) - 3 \min(b_4, c_4) - \sum_{i \geqslant 5} (i-1) (b_i + c_i) - 1.
\end{equation}
Then by (80) we have
\begin{equation}
\nu \leqslant \max(\nu_1, \nu_2) + \epsilon.
\end{equation}
By (18), (23), (67) and (78), we know that
\begin{align}
\nonumber \nu_1 =&\ a_3 + b_3 + \min(b_4, c_4) + b_5 + c_5 + \sum_{i \geqslant 6}(b_i + c_i) \\
\nonumber \leqslant&\ a_3 + b_3 + \min(b_4, c_4) + b_5 + c_5 + \frac{1}{5}\sum_{i \geqslant 6}(i - 1)(b_i + c_i) \\
\nonumber \leqslant&\ a_3 + b_3 + \frac{b_4 + c_4}{2} + b_5 + c_5 + \frac{1}{5}\left(\frac{4}{3} +\delta_{bc} - (b_2 + c_2) - 2 (b_3 + c_3) - 3 (b_4 + c_4) - 4 (b_5 + c_5) \right) \\
\nonumber \leqslant&\ a_3 + b_3 + \frac{b_4 + c_4}{2} + b_5 + c_5 + \frac{1}{5}\left(\frac{4}{3} +\delta_{bc} - 3 (b_4 + c_4) - 4 (b_5 + c_5) \right) \\
\nonumber \leqslant&\ a_3 + b_3 + \frac{1}{5}\left(\frac{4}{3} +\delta_{bc} \right) - \frac{b_4 + c_4}{10} + \frac{b_5 + c_5}{5} \\
\nonumber \leqslant&\ a_3 + b_3 + \frac{1}{5}\left(\frac{4}{3} +\delta_{bc} + (b_5 + c_5) \right) \\
\nonumber \leqslant&\ (1 - \theta) + \frac{1}{5}\left(\frac{4}{3} + \left(\frac{2}{3} -\theta \right) + (6 - 9\theta) \right) \\
=&\ \frac{13}{5} - 3 \theta < \theta - \epsilon
\end{align}
when $\theta > \frac{13}{20} = 0.65$.

Note that we have
\begin{equation}
\sum_{i} i s_i - 1 \leqslant \sum_{i} i s_i - \sum_{i} i a_i =  \sum_{i} i (b_i + c_i)
\end{equation}
by (8). Then by (16), (19), (75), (78), (85) and assumptions, for $\nu_2$ we have
\begin{align}
\nonumber \nu_2 =&\ \left( \sum_{i} i s_i - 1\right) - 2(a_3 + b_3) - 3 \min(b_4, c_4) - \sum_{i \geqslant 5} (i-1) (b_i + c_i) \\
\nonumber \leqslant&\ \left( \sum_{i} i (b_i + c_i) - \sum_{i \geqslant 5} (i-1) (b_i + c_i) \right) - 2(a_3 + b_3) - 3 \min(b_4, c_4) \\
\nonumber =&\ \left( \sum_{i} (b_i + c_i) + \sum_{i \leqslant 4} (i-1) (b_i + c_i) \right) - 2(a_3 + b_3) - 3 \min(b_4, c_4) \\
\nonumber =&\ \sum_{i} (b_i + c_i) + (b_2 + c_2) + 2 (b_3 + c_3) + 3 (b_4 + c_4) - 2(a_3 + b_3) - 3 \min(b_4, c_4) \\
\nonumber =&\ \sum_{i} (b_i + c_i) + (b_2 + c_2) - 2(a_3 - c_3) + 3 \max(b_4, c_4) \\
\nonumber \leqslant&\ \frac{2}{3} - \delta_{bc} + s_2 - 2 a_3 + 2 c_3 + 3 \max(b_4, c_4) \\
\nonumber \leqslant&\ \frac{2}{3} - 2 \left(\theta -\frac{2}{3} \right) + \left( \frac{49}{12} - \frac{23}{4}\theta \right) - 2 (2 \theta - 1) + 2 (2 - 3 \theta) + 3 \left(3 - \frac{9}{2} \theta \right) \\
=&\ \frac{253}{12} - \frac{125}{4} \theta < \theta - \epsilon
\end{align}
when $\theta > \frac{253}{387} \approx 0.6537$. Now by (83), (84) and (86), we get the desired result.

\subsubsection{Subcase 2.2: $b_3 + c_3 < \frac{1}{2} \theta - \frac{1}{2} s_2$} Now by (31) and a similar discussion as in (44)--(45), we can assume that
\begin{equation}
b_3 + c_3 \leqslant 1 - \theta - s_1 - s_2.
\end{equation}
By (26), (37) and (87) we know that
\begin{align}
\nonumber a_3 = s_3 - (b_3 + c_3) \geqslant&\ \left(1 - 2 s_1 - \frac{3}{2} s_2 - \frac{1}{2} s_4 - \frac{5}{2} \delta_s \right) - (1 - \theta - s_1 - s_2) \\
\nonumber =&\ \theta - \frac{5}{2} \delta_s - s_1 - \frac{1}{2}(s_2 + s_4) \\
\nonumber \geqslant&\ \theta - \frac{5}{2} \left(1-\frac{3}{2}\theta \right) - \frac{1}{2}\left(\frac{3}{2}-\frac{3}{2}\theta \right) - s_1 \\
=&\ \left(\frac{11}{2}\theta -\frac{13}{4} \right) - s_1.
\end{align}
Note that
\begin{equation}
\frac{11}{2}\theta -\frac{13}{4} > 1 - \theta
\end{equation}
since $\theta > \frac{17}{26} \approx 0.6538$, we have
\begin{equation}
a_3 > 1 - \theta - s_1.
\end{equation}
By (31) and a similar discussion as in (44)--(45), we can assume that
\begin{equation}
a_3 > \frac{1}{2}\theta.
\end{equation}
Note that
\begin{equation}
\frac{1}{2}\theta > 2 \theta - 1
\end{equation}
since $\theta < \frac{2}{3}$, we have
\begin{equation}
a_3 > 2 \theta - 1.
\end{equation}
Now by the discussions in \textit{Subcase 2.1}, we get the desired result.

Now, we will prove \textbf{Case 2} by showing that (12) holds for any $(s_1, s_2) \in [0,1]^2$ (with the assumption $s_2 < k$). We shall consider the following 4 subcases:
\begin{align}
\begin{cases}
(2.3) & 4 s_1 + 3 s_2 > 4 - 5 \theta, \\
(2.4) & 4 s_1 + s_2 < 37 \theta - 24, \\
(2.5) & 6 - 9 \theta \leqslant s_2 \leqslant \frac{7}{3} \theta - \frac{4}{3}, \\
(2.6) & 2 s_1 - s_2 > 2 - 3 \theta.
\end{cases}
\end{align}
Note that every point in $[0,1]^2$ is covered by one of the above 4 subcases when $\theta \geqslant \frac{23}{35} \approx 0.6571$. If $\theta < \frac{23}{35}$, there are two triangles that are not covered by any of the cases.

\subsubsection{Subcase 2.3: $4 s_1 + 3 s_2 > 4 - 5 \theta$}
By (73) we know that
\begin{equation}
b_3 + c_3 \leqslant 2 - 2\theta - 2 s_1 - 2 s_2.
\end{equation}
By the assumption we know that
\begin{equation}
- 2 s_1 - 2 s_2 < \frac{5}{2} \theta - 2 - \frac{1}{2} s_2
\end{equation}
Now, by (95) and (96) we have
\begin{equation}
b_3 + c_3 < 2 - 2\theta + \frac{5}{2} \theta - 2 - \frac{1}{2} s_2 = \frac{1}{2} \theta - \frac{1}{2} s_2.
\end{equation}
Hence \textit{Subcase 2.2} completes the proof.

\subsubsection{Subcase 2.4: $4 s_1 + s_2 < 37 \theta - 24$}
By (24) and (74) we know that
\begin{align}
\nonumber b_3,\ c_3 \leqslant&\ 1 - \theta - a_3 \\
\nonumber \leqslant&\ 1 - \theta - (2 \theta - 1 - s_1 - 4 \delta_s) \\
\nonumber =&\ 2 - 3\theta + s_1 + 4 \delta_s \\
\nonumber =&\ 2 - 3\theta + s_1 + 4 \left(1 -\frac{3}{2} \theta \right) \\
=&\ 6 - 9 \theta + s_1, \\
b_3 + c_3 \leqslant&\ 12 - 18 \theta + 2 s_1.
\end{align}
By the assumption we know that
\begin{equation}
2 s_1 < \frac{37}{2} \theta - 12 - \frac{1}{2} s_2.
\end{equation}
Now, by (99) and (100) we have
\begin{equation}
b_3 + c_3 \leqslant 12 - 18 \theta + 2 s_1 < 12 - 18 \theta + \frac{37}{2} \theta - 12 - \frac{1}{2} s_2 = \frac{1}{2} \theta - \frac{1}{2} s_2.
\end{equation}
Hence \textit{Subcase 2.2} completes the proof.

\subsubsection{Subcase 2.5: $6 - 9 \theta \leqslant s_2 \leqslant \frac{7}{3} \theta - \frac{4}{3}$}
By (22) and (74) we know that
\begin{equation}
a_3 \geqslant 2 \theta - 1 - s_1 - 4 \delta_s \geqslant 2 \theta - 1 - s_1 - 4 \left(1 - \frac{3}{2}\theta \right) \geqslant 8 \theta - 5 - s_1.
\end{equation}
If $6 - 9 \theta \leqslant s_2$, we have
\begin{equation}
a_3 \geqslant 8 \theta - 5 - s_1 \geqslant 8 \theta - 5 - s_1 + (6 - 9 \theta - s_2) = 1 - \theta - s_1 - s_2.
\end{equation}
By (31) and a similar discussion as in (44)--(45) we can assume
\begin{equation}
a_3 \geqslant \frac{1}{2} \theta - \frac{1}{2} s_2.
\end{equation}
Now, (24) yields
\begin{align}
\nonumber b_3,\ c_3 \leqslant&\ 1 - \theta - a_3 \\
\nonumber \leqslant&\ 1 - \theta - \left( \frac{1}{2} \theta - \frac{1}{2} s_2 \right) \\
=&\ 1 - \frac{3}{2} \theta + \frac{1}{2} s_2, \\
b_3 + c_3 \leqslant&\ 2 - 3 \theta + s_2.
\end{align}
If $s_2 \leqslant \frac{7}{3} \theta - \frac{4}{3}$, we have
\begin{equation}
s_2 \leqslant \frac{7}{2} \theta - 2 - \frac{1}{2} s_2.
\end{equation}
Now, by (106) and (107) we have
\begin{equation}
b_3 + c_3 \leqslant 2 - 3 \theta + s_2 \leqslant 2 - 3 \theta + \frac{7}{2} \theta - 2 - \frac{1}{2} s_2 = \frac{1}{2} \theta - \frac{1}{2} s_2.
\end{equation}
Hence \textit{Subcase 2.2} completes the proof.

\subsubsection{Subcase 2.6: $2 s_1 - s_2 > 2 - 3 \theta$}
In this case the two intervals in (31) overlap, hence we have $\nu \leqslant \theta$ if
\begin{equation}
\tau_3 \in \left(1-\theta -s_1 -s_2,\ \frac{1}{2}\theta  \right).
\end{equation}
In \textit{Subcases 2.3 and 2.4} we prove the cases $4 s_1 + 3 s_2 > 4 - 5 \theta$ and $4 s_1 + s_2 < 37 \theta - 24$, so we can assume that
\begin{equation}
4 s_1 + 3 s_2 \leqslant 4 - 5 \theta 
\end{equation}
and
\begin{equation}
4 s_1 + s_2 \geqslant 37 \theta - 24.
\end{equation}
By (110) we have
\begin{equation}
s_1 \leqslant \frac{4 - 5\theta}{4} = 1 - \frac{5}{4} \theta
\end{equation}
and
\begin{equation}
s_2 \leqslant \frac{1}{3} (4 - 5 \theta - 4 s_1).
\end{equation}
Now, (111) and (113) give that
\begin{align}
\nonumber s_2 \leqslant&\ \frac{1}{3} (4 - 5 \theta - 4 s_1) \\
\nonumber \leqslant&\ \frac{1}{3} (4 - 5 \theta - (37 \theta - 24 - s_2) ) \\
=&\ \frac{1}{3} (28 - 42 \theta + s_2 ), \\
s_2 \leqslant&\ 14 - 21 \theta.
\end{align}
Note that
\begin{equation}
14 - 21 \theta < \frac{7}{3} \theta - \frac{4}{3}
\end{equation}
when $\theta > \frac{23}{35} \approx 0.6571$, we have
\begin{equation}
s_2 \leqslant \frac{7}{3} \theta - \frac{4}{3}.
\end{equation}
If $6 - 9 \theta \leqslant s_2 \leqslant \frac{7}{3} \theta - \frac{4}{3}$, by \textit{Subcase 2.5} we have the desired result. Otherwise we have
\begin{equation}
s_2 \leqslant 6 - 9 \theta.
\end{equation}
By the result proved in \textit{Subcase 2.1}, we can also assume that $a_3 < 2\theta - 1$. Since $2 \theta - 1 < \frac{1}{2} \theta$ when $\theta < \frac{2}{3}$, by (31) and a similar discussion as in (44)--(45), we have
\begin{equation}
a_3 < 1-\theta -s_1 -s_2.
\end{equation}
We shall consider the following two cases.

\textit{Subcase 2.6.1: $b_3 + c_3 < 1-\theta -s_1 -s_2$.} In this case we have, by the assumption and (119),
\begin{equation}
s_3 < 2- 2\theta - 2 s_1 - 2 s_2.
\end{equation}
Now, by (36) and (120) we have
\begin{equation}
1 - 3 s_1 - 2 s_2 - 4 \delta_s < 2- 2\theta - 2 s_1 - 2 s_2
\end{equation}
and thus
\begin{equation}
2 \theta - 1 - 4 \delta_s < s_1.
\end{equation}
By (112) and (122), we have
\begin{equation}
2 \theta - 1 + 4 \delta < 1 - \frac{5}{4} \theta,
\end{equation}
which holds true only when $\theta < \frac{8}{13} \approx 0.6154$. This contradicts with our value of $\theta$.

\textit{Subcase 2.6.2: $b_3 + c_3 \geqslant 1-\theta -s_1 -s_2$.} By the assumption and (109), after a similar discussion as in (44)--(45) we have
\begin{equation}
b_3 + c_3 \geqslant \frac{1}{2}\theta.
\end{equation}
Since $b_3 > c_3$, we have $b_3 > \frac{1}{4} \theta$. Now $a_3 > b_3$ yields
\begin{equation}
\frac{1}{4} \theta < b_3 < a_3 < 1-\theta -s_1 -s_2,
\end{equation}
\begin{equation}
s_1 + s_2 < 1 - \frac{5}{4} \theta.
\end{equation}

By Lemma~\ref{Fourier}, we know that
\begin{equation}
\nu \leqslant \frac{1}{2}\left(1 + \delta + \sum_{i \neq 3}\max(a_i, b_i) \right).
\end{equation}
Note that
\begin{equation}
\sum_{i} (\max(a_i, b_i) + \min(a_i, b_i)) = \sum_{i}(a_i + b_i) = \frac{2}{3} - \delta_{ab},
\end{equation}
we have
\begin{align}
\nonumber 2 \nu - 1 - \delta \leqslant&\ \sum_{i \neq 3}\max(a_i, b_i) \\
\nonumber \leqslant&\ \max(a_1, b_1) + \max(a_2, b_2) + \sum_{i \geqslant 4}(\max(a_i, b_i) + \min(a_i, b_i)) \\
\nonumber =&\ \max(a_1, b_1) + \max(a_2, b_2) + \left(\frac{2}{3} - \delta_{ab} - \sum_{i \leqslant 3}(\max(a_i, b_i) + \min(a_i, b_i)) \right) \\
=&\ \frac{2}{3} - \delta_{ab} - \min(a_1, b_1) - \min(a_2, b_2) - \min(a_3, b_3) - \max(a_3, b_3).
\end{align}
By (34) we know that
\begin{equation}
3 a_1 + 2 a_2 + a_3 \geqslant \frac{1}{3} - 4 \delta_a
\end{equation}
and
\begin{equation}
3 b_1 + 2 b_2 + b_3 \geqslant \frac{1}{3} - 4 \delta_b.
\end{equation}
Thus,
\begin{equation}
3 \min(a_1, b_1) \geqslant \frac{1}{3} - 2 \max(a_2, b_2) - \max(a_3, b_3) - 4 \max(\delta_{a}, \delta_{b}).
\end{equation}
Now, by (19), (125) and (132) we have
\begin{align}
\nonumber 2 \nu - 1 - \delta \leqslant&\ \frac{2}{3} - \delta_{ab} - \min(a_1, b_1) - \min(a_2, b_2) - \min(a_3, b_3) - \max(a_3, b_3) \\
\nonumber \leqslant&\ \frac{2}{3} - \delta_{ab} - \frac{1}{3}\left(\frac{1}{3} - 2 \max(a_2, b_2) - \max(a_3, b_3) - 4 \max(\delta_{a}, \delta_{b}) \right) \\
\nonumber & - \min(a_2, b_2) - \min(a_3, b_3) - \max(a_3, b_3) \\
\nonumber \leqslant&\ \frac{5}{9} + \frac{2}{3}\max(a_2, b_2) + \left(\frac{4}{3}\max(\delta_{a}, \delta_{b}) - \delta_{ab} \right) -\left(\min(a_3, b_3) + \frac{2}{3}\max(a_3, b_3) \right) \\
\nonumber \leqslant&\ \frac{5}{9} + \frac{2}{3}\max(a_2, b_2) + \frac{1}{3}\max(\delta_{a}, \delta_{b}) - \min(\delta_{a}, \delta_{b}) -\frac{5}{3} \min(a_3, b_3) \\
\nonumber \leqslant&\ \frac{5}{9} + \frac{2}{3}\max(a_2, b_2) + \frac{1}{3} \left(\frac{4}{3}-2\theta \right) - \left(\theta-\frac{2}{3} \right) -\frac{5}{3} \left(\frac{1}{4}\theta \right) \\
=&\ \frac{5}{3} - \frac{25}{12}\theta + \frac{2}{3}\max(a_2, b_2), \\
\nu \leqslant &\ \frac{4}{3} - \frac{25}{24}\theta + \frac{1}{3}\max(a_2, b_2) + \frac{1}{2}\delta.
\end{align}
By (134), we know that (12) holds if we have
\begin{equation}
\max(a_2, b_2) < \frac{49}{8} \theta - 4.
\end{equation}
Now we assume that
\begin{align}
\max(a_2, b_2) \geqslant \frac{49}{8} \theta - 4.
\end{align}
By similar arguments as above, we also have
\begin{align}
\max(a_2, c_2) \geqslant \frac{49}{8} \theta - 4
\end{align}
and
\begin{align}
\max(b_2, c_2) \geqslant \frac{49}{8} \theta - 4,
\end{align}
which mean that at least two of $a_2, b_2, c_2$ are $\geqslant \frac{49}{8} \theta - 4$, but then we have
\begin{equation}
s_2 \geqslant 2\left(\frac{49}{8} \theta - 4 \right) = \frac{49}{4} \theta - 8,
\end{equation}
which is larger than $6 - 9 \theta$ when $\theta > \frac{56}{85}$ and thus contradicts with (118). That is why we stop at this point.

Finally, combining all above cases, Theorem~\ref{t3} is proved.

\bibliographystyle{plain}
\bibliography{bib}

\end{document}